## RESEARCH

# Transcendence of certain sequences of algebraic numbers

Mathias L. Laursen

*Correspondence:
mll@math.au.dk
Department of Mathematics,
Aarhus University, Ny
Munkegade 118, 8000 Aarhus C,
Denmark

**Abstract**

Using Schmidt's Subspace Theorem, this paper improves and extends an existing transcendence result for sequences of algebraic numbers. The theorems thus produced correspond to a central theorem on the irrationality of sequences due to Erdős.

## 1 Introduction and main results

Proving whether a given real number is algebraic, or even rational, can be a quite frustrating endeavour. While more than a century and a half have passed since Hermite proved that the number $e = \sum_{n=0}^{\infty} \frac{1}{n!}$ is transcendental in early 1873, it remains unsolved if the number $\sum_{n=0}^{\infty} \frac{1}{n!+1}$ is even irrational, despite what may appear as a much similar construction. Similarly, it is well-known that the Riemann $\zeta$ function defined as $\zeta(s) = \sum_{n=1}^{\infty} n^{-s}$ for $\Re(s) > 1$ is transcendental when $s$ is a positive even integer while the question of irrationality remains open when $s \geq 5$ is any fixed odd integer. In other words, we have a multitude of interesting numbers that we know to be transcendental but where a small perturbation to the infinite series used to describe them renders even the question of irrationality exceedingly hard to settle. Aiming away from frustrations of this kind, this paper studies irrationality and transcendence criteria that are less sensitive to such perturbations.

Following the notions of Erdős and Graham [1] (respectively Hančl [2]), we say that a sequence $\{a_n\}_{n=1}^{\infty}$ of real or complex numbers is irrational (respectively transcendental) if the sum of the series $\sum_{n=1}^{\infty} \frac{1}{a_n c_n}$ is irrational (respectively transcendental) for any sequence $\{c_n\}_{n=1}^{\infty}$ of positive integers. An early and central result on irrational sequences was proven in 1975 by Erdős [3].

**Theorem 1.1** (Erdős) *Let $\varepsilon > 0$, and let $\{a_n\}_{n=1}^{\infty}$ be an increasing sequence of integers such that $a_n \geq n^{1+\varepsilon}$ for all n. Suppose*

$$\limsup_{n \to \infty} a_n^{2^{-n}} = \infty.$$

*Then the sequence $\{a_n\}_{n=1}^\infty$ is irrational.*

As noted in [3], the number 2 in the theorem is best possible, in the sense that there exist sequences $\{a_n\}_{n=1}^\infty$ of positive integers that satisfy $\limsup_{n\to\infty} a_n^{A^{-n}} = \infty$ for all $A < 2$ while the sum of the series $\sum_{n=1}^\infty 1/a_n$ is rational. Still, much effort has been applied to extend this result (see [4] for a broader overview). One such result is the below theorem by Hančl [5], which gives a corresponding condition for a sequence of (not necessarily integral) rational numbers to be transcendental.

**Theorem 1.2** (Hančl) *Let $\gamma > 2\varepsilon > 0$ and $1 > \alpha > \frac{\log(3+2\varepsilon)}{\log(3+\gamma)}$, and let $\{a_n\}_{n=1}^\infty$ and $\{b_n\}_{n=1}^\infty$ be sequences of positive integers such that*

$$n^{1+\varepsilon} < a_n \leq a_{n+1}, \quad \limsup_{n\to\infty} a_n^{(3+\gamma)^{-n}} = \infty,$$

*and*

$$b_n < a_n^{\varepsilon/(1+\varepsilon)} 2^{-\log_2^\alpha a_n}. \tag{1}$$

*Then the sequence $\{a_n/b_n\}_{n=1}^\infty$ is transcendental.*

One way to extend on this result is to broaden the family of numbers that may be contained by the sequence to a greater class of algebraic numbers. The best known result in this direction is due to Andersen and Kristensen [6], which gives sufficient conditions for bounding the algebraic degree from below. In their theorem, they use the notion of algebraic integers, which are defined as the algebraic numbers that have a monic minimal polynomial over the integers. Given an algebraic extension $K$ of $\mathbb{Q}$, we will use $\mathcal{O}_K$ to denote the set of algebraic integers contained in $K$. Recall that $\mathcal{O}_K$ forms a subring of $K$. They also use the notion of a *house*, written as $\overline{|a|} := \max_{1 \leq j \leq d} |a^{(j)}|$, where $a^{(1)}, \ldots, a^{(d)}$ denotes the conjugates of an algebraic number $a$, i.e., the roots of its minimal polynomial over $\mathbb{Z}$. In terms of irrational and transcendental sequences, the theorem reduces to the below result.

**Theorem 1.3** (Andersen–Kristensen) *Let $d \in \mathbb{N}$ be a positive integer, and let $\{a_n\}_{n=1}^\infty$ be a sequence of algebraic integers of degree $\deg a_n \leq d$ such that*

$$n^{1+\varepsilon} \leq \overline{|a_n|} = |a_n| \leq |a_{n+1}|.$$

*Suppose $\Re(a_n) > 0$ for all $n$ or that $\Im(a_n) > 0$ for all $n$. If*

$$\limsup_{n\to\infty} |a_n|^{\prod_{i=1}^{n-1}(d^i+d)^{-1}} = \infty,$$

*then $\{a_n\}_{n=1}^\infty$ is irrational. Furthermore, if for all $D \in \mathbb{N}$,*

$$\limsup_{n\to\infty} |a_n|^{D^{-n}\prod_{i=1}^{n-1}(d^i+d)^{-1}} = \infty,$$

*then $\{a_n\}_{n=1}^\infty$ is transcendental.*

Note that the restriction on $\Re(a_n)$ or $\Im(a_n)$ corresponds to the restriction that $a_n$ be positive in Theorem 1.1. Furthermore, as can be seen in the proof, the somewhat extensive assumptions on the divergence of the limsup of $a_n$ is in part due to the fact that each successive $a_N$ may potentially increase the algebraic degree of $\sum_{n=1}^N 1/a_n$ by a factor of $d$. By assuming the $a_n$ to come from a fixed number field $K$, the limsup conditions would thus be much weakened, replacing the product $\prod_{i=1}^{n-1}(d^i + d)^{-1}$ with $(2d)^{-n}$. The main results of this paper are improvements to this result when all $a_n$ are contained in the same number field, in terms of allowing sequences with non-integral elements, replacing the restriction $\overline{a_n} = |a_n|$ with much weaker conditions, and weakening the limsup criteria. In our first result, below, we assume $a_n$ to be rational but allow $b_n$ to attain certain algebraic irrational values. This will be proven in Sect. 4.

**Theorem 1.4** *Let $K$ be a number field of dimension $d \geq 2$, and let $x_1, \ldots, x_D \in K$. Consider real numbers $\varepsilon > 0$, $0 < \alpha < 1 \leq y$, and $\beta \in [0, \frac{\varepsilon}{1+\varepsilon})$. Let $\{a_n\}_{n=1}^{\infty}$ be a sequence of positive integers such that*

$$n^{1+\varepsilon} \leq a_n \leq a_{n+1}, \tag{2}$$

*and let $\{b_n\}_{n=1}^{\infty}$ be a sequence of non-zero numbers so that $b_n = \sum_{i=1}^{D} b_{i,n} x_i$ for suitable $b_{i,n} \in \mathbb{Z}$. Suppose for $n$ sufficiently large, $i = 1, \ldots, D$, and some fixed $\zeta \in \mathbb{C}$,*

$$|b_n| \leq a_n^{\beta} 2^{\log_2^{\alpha} a_n}, \quad |b_{i,n}| \leq a_n^{y} 2^{\log_2^{\alpha} a_n}, \tag{3}$$

*and*

$$\Re(\zeta b_n) > 0. \tag{4}$$

*Then the sequence $\{a_n/b_n\}_{n=1}^{\infty}$ is irrational if*

$$\limsup_{n \to \infty} a_n^{\left(\frac{dy}{1-\beta}+1\right)^{-n}} = \infty,$$

*and it is transcendental if*

$$\limsup_{n \to \infty} a_n^{\left(\frac{d^2 y}{1-\beta}+1\right)^{-n}} = \infty. \tag{5}$$

*Remark 1.5* In the proof of the theorem, assumption (4) will only be used once and only to ensure that the partial sums $\sum_{n=1}^{N} \frac{b_n}{a_n}$ do not take the same value infinitely often (see the proof of Proposition 4.3 later in this paper). Therefore, assumption (4) can be replaced by any other assumption that preserves this property.

The main novelty of this result is the improved transcendence criterion, which relies on Schmidt's Subspace Theorem along with ideas from [7] to exclude near all algebraic numbers as possible values for the sum of $\sum_{n=1}^{\infty} \frac{1}{a_n c_n}$, leaving only a finite field of potential algebraic values to be dealt with in the spirit of [6].

Since the above theorem assumes $a_n$ to be rational, much of the arithmetic information regarding the number $a_n/b_n$ – such as its algebraic degree – is carried solely by $b_n$. This is in some contrast to Theorems 1.1 and 1.3, which can be viewed as having $b_n$ constantly 1, so that all arithmetic information is stored in $a_n$ alone. By modifying the proof of Theorem 1.4 in order to get a result where we again have most of the arithmetic information carried by $a_n$, we reach the below result. Unfortunately, this version of the theorem is a bit more complicated to read, which is in part due to the method of proof as it requires $\sum_{n=1}^{N} \frac{b_n}{a_n c_n}$ to be written as a $\mathbb{Q}$-linear combination of the $x_i$ – something that is more easily and neatly done when the $a_n$ are guaranteed to be rational.

In the theorem and for the rest of this paper, $\mathcal{N} : \overline{\mathbb{Q}} \to \mathbb{Q}$ denotes the map that sends each algebraic number to the product of its algebraic conjugates, and $\mathcal{N}_K : K \to \mathbb{Q}$ denotes the field norm for the finite extension $K \subseteq \mathbb{Q}$. Notice that $\mathcal{N}_K(a) = \mathcal{N}(a)^{d/\deg a}$ for all $a \in K$, where $d$ denotes the degree of the extension $K \subseteq \mathbb{Q}$.

**Theorem 1.6** *Let $K$ be a number field of dimension $d$, and let $x_1, \ldots, x_D \in K$. Consider real numbers $\alpha, \delta, \varepsilon > 0$, $\beta, \eta_1 \geq 0$, and $\eta_2, y \geq 1$ such that $\alpha < 1$, $\beta < \varepsilon/(1+\varepsilon)$, and $\eta_1 \leq (d-1)y + \beta$. Let $\{a_n\}_{n=1}^\infty$ be a sequence of non-zero numbers given by $a_n = \sum_{i=1}^D a_{i,n} x_i$ with $a_{i,n} \in \mathbb{Z}$ such that for all sufficiently large $n$,*

$$n^{1+\varepsilon} \leq |a_n| \leq |a_{n+1}|, \tag{6}$$
$$|\mathcal{N}_K(a_n)| \geq |a_n|^{\eta_1} 2^{-\log_2^\alpha |a_n|}, \tag{7}$$

*and*

$$r_n |\mathcal{N}(a_n/r_n)| \leq |a_n|^{\eta_2} 2^{\log_2^\alpha |a_n|}, \tag{8}$$

*where $r_n := \gcd(a_{1,n}, \ldots, a_{D,n})$. Let $\{b_n\}_{n=1}^\infty$ be a sequence of positive integers such that for some fixed $\zeta \in \mathbb{C}$, each $i = 1, \ldots, D$, and all sufficiently large $n$,*

$$b_n \leq |a_n|^\beta 2^{\log_2^\alpha |a_n|}, \quad |a_{i,n}| \leq |a_n|^y 2^{\log_2^\alpha |a_n|},$$

*and $\Re(\zeta a_n) > 0$. Then the sequence $\{a_n/b_n\}_{n=1}^\infty$ is irrational if*

$$\limsup_{n \to \infty} |a_n|^{\left(\frac{d(y+\beta)}{1-\beta}+1\right)^{-n}} = \infty,$$

*and it is transcendental if*

$$\limsup_{n \to \infty} |a_n|^{\left(\frac{\eta_2 + d((d-1)y+\beta+\eta_2-\eta_1)+\delta}{1-\beta}+1\right)^{-n}} = \limsup_{n \to \infty} |a_n|^{\left(\frac{d^2(y+\beta)}{1-\beta}+1\right)^{-n}} = \infty.$$

In the proof of Theorem 1.6, it makes little difference if we also allow $b_n$ to be irrational. Doing so leads to the below generalization, which we will prove in Sect. 5.

**Theorem 1.7** *Let $K$ be a number field of dimension $d$, and let $x_1, \ldots, x_D \in K$. Consider real numbers $\alpha, \delta, \varepsilon > 0$, $\beta, \eta_1 \geq 0$, $\eta_2, y_1 \geq 1$, and $y_2 \geq \beta$ such that $\alpha < 1$, $\beta < \varepsilon/(1+\varepsilon)$, and $\eta_1 \leq (d-1)y_1 + y_2$. Let $\{a_n\}_{n=1}^\infty$ and $\{b_n\}_{n=1}^\infty$ be sequences of non-zero numbers given by $a_n = \sum_{i=1}^D a_{i,n} x_i$ and $b_n = \sum_{i=1}^D b_{i,n} x_i$ with $a_{i,n}, b_{i,n} \in \mathbb{Z}$ such that inequalities (6), (7), and (8) are satisfied for $n$ sufficiently large. For each $i = 1, \ldots, D$, and some fixed $\zeta \in \mathbb{C}$, suppose additionally that*

$$|a_{i,n}| \leq |a_n|^{y_1} 2^{\log_2^\alpha |a_n|}, \quad |b_{i,n}| \leq |a_n|^{y_2} 2^{\log_2^\alpha |a_n|}, \tag{9}$$
$$|b_n| \leq |a_n|^\beta 2^{\log_2^\alpha |a_n|}, \tag{10}$$

*and*

$$\Re(\zeta a_n/b_n) > 0, \tag{11}$$

*when $n$ is sufficiently large. Then the sequence $\{a_n/b_n\}_{n=1}^\infty$ is irrational if*

$$\limsup_{n \to \infty} |a_n|^{\left(\frac{d(y_1+y_2)}{1-\beta}+1\right)^{-n}} = \infty,$$

*and it is transcendental if*

$$\limsup_{n \to \infty} |a_n|^{\left(\frac{\eta_2 + d((d-1)y_1 + y_2 + \eta_2 - \eta_1) + \delta}{1-\beta} + 1\right)^{-n}} = \limsup_{n \to \infty} |a_n|^{\left(\frac{d^2(y_1 + y_2)}{1-\beta} + 1\right)^{-n}} = \infty. \quad (12)$$

*Remark 1.8* Similarly to Theorem 1.4, the assumption (11) can be replaced with any other assumption that ensures that the partial sums $\sum_{n=1}^{N} \frac{b_n}{a_n}$ do not have the same value infinitely often (see the proof of Proposition 5.3 later in this paper).

*Remark 1.9* Suppose that $\{a_n\}_{n=1}^{\infty}$ and $\{b_n\}_{n=1}^{\infty}$ satisfy the assumptions of either of Theorems 1.4, 1.6, and 1.7 for some choice of $x_1, \ldots, x_D$. If $x_1', \ldots, x_{D'}' \in K$ such that $a_n$ and $b_n$ lie in the $\mathbb{Q}$-linear span of these numbers, then $\{a_n\}_{n=1}^{\infty}$ and $\{b_n\}_{n=1}^{\infty}$ satisfy the assumptions of the same theorem with $x_1'/Q, \ldots, x_{D'}'/Q$ instead of $x_1, \ldots, x_D$, where $Q$ is a positive integer that depends only on $x_1, \ldots, x_D$ and $x_1', \ldots, x_{D'}'$.

To see that this is indeed the case, pick one of Theorems 1.4, 1.6, and 1.7, use the notation from that theorem, and assume the conditions are satisfied. Let $d'$ denote the dimension of the $\mathbb{Q}$-linear span of $x_1', \ldots, x_{D'}'$. By renumbering if necessary, we may assume that $x_1', \ldots, x_{d'}'$ are linearly independent. Let $\tilde{x}_1, \ldots \tilde{x}_d$ be a $\mathbb{Q}$-linear basis of $K$ so that $\tilde{x}_j = x_j'$ for $1 \leq j \leq d'$, and write $x_i = \sum_{j=1}^{d} \frac{p_{i,j}}{q_{i,j}} \tilde{x}_j$, for suitable choices of $p_{i,j} \in \mathbb{Z}$ and $q_{i,j} \in \mathbb{N}$. Pick $Q = \prod_{i=1}^{D} \prod_{j=1}^{d} q_{i,j}$. Let $\xi$ denote either letter $a$ or $b$ so that $\xi_n$ is not assumed to be a positive integer by the chosen theorem. Seeing that

$$\xi_n = \sum_{i=1}^{D} \xi_{i,n} x_i = \sum_{j=1}^{d} \sum_{i=1}^{D} \frac{Q}{q_{i,j}} p_{i,j} a_{i,n} \frac{\tilde{x}_j}{Q},$$

write $\tilde{\xi}_{j,n} = \sum_{i=1}^{D} \frac{Q}{q_{i,j}} p_{i,j} \xi_{i,n}$, and set $\xi_{j,n}' = \tilde{\xi}_{j,n}$ for $1 \leq j \leq d'$ and $\xi_{j,n}' = 0$ for $d' < j \leq D'$. Since each $\xi_n$ is contained in the span of $x_1', \ldots, x_{D'}'$ (and so in the span of $\tilde{x}_1 = x_1', \ldots, \tilde{x}_{d'} = x_{d'}'$), it follows that $\tilde{\xi}_{j,n} = 0$ for $d' < j \leq d$, and so

$$\xi_n = \sum_{j=1}^{d} \tilde{\xi}_{j,n} \frac{\tilde{x}_j}{Q} = \sum_{j=1}^{D'} \xi_{j,n}' \frac{x_j'}{Q},$$

while

$$\max_{1 \leq j \leq D'} |\xi_{j,n}'| = \max_{1 \leq j \leq d} |\tilde{\xi}_{j,n}| \leq \sum_{i=1}^{D} \frac{Q}{q_{i,j}} |p_{i,j}| |\xi_{i,n}| \leq DQ \max_{i,j} |p_{i,j}| \max_{1 \leq i \leq D} |\xi_{i,n}|.$$

By replacing $\alpha$ with $(1 + \alpha)/2$, it follows that $\{a_n\}_{n=1}^{\infty}$ and $\{b_n\}_{n=1}^{\infty}$ satisfy the assumptions of the chosen theorem with $x_1'/Q, \ldots, x_{D'}'/Q$ instead of $x_1, \ldots, x_D$.

Notice that the sum $y_1 + y_2$ in Theorem 1.7 corresponds to $y + \beta$ in Theorem 1.4. As such, one should not expect to be able to derive Theorem 1.4 as a corollary to Theorem 1.7. This is further underlined by Example 2.3 in section 2. Similarly, as will be seen from Example 2.7, there are cases where Theorem 1.7 is applicable while the other two theorems are not.

## 2 Examples

We will now go through a few applications of the main theorems in order to better understand the strengths and differences of applicability between them. For this purpose, we will say that a theorem is immediately applicable to a sequence $\{x_n\}_{n=1}^\infty$ if there are sequences $\{a_n\}_{n=1}^\infty$ and $\{b_n\}_{n=1}^\infty$ that satisfy $x_n = a_n/b_n$ and the assumptions of the theorem.

For these examples, we make use of the Fibonacci sequence $F_n$, defined by $F_0 = 0$, $F_1 = 1$, and $F_{n+1} = F_n + F_{n-1}$, along with the golden ratio $\varphi = (1 + \sqrt{5})/2$ and its conjugate $\bar{\varphi} = (1 - \sqrt{5})/2 = -\varphi^{-1}$. Recall that $\varphi^n = F_n\varphi + F_{n-1}$ and $\bar{\varphi}^n = F_n\bar{\varphi} + F_{n-1}$ for each $n \in \mathbb{N}$.

The first example, below, shows the strengths of Theorem 1.4 in terms of providing transcendence of $K$-linear combinations of multiple series of rational numbers when $K$ is a suitable number field.

*Example 2.1* Let $x$ be any algebraic number of degree at most 2, and let $\{c_n\}_{n=1}^\infty$ be a sequence of positive integers. Then

$$x \sum_{n=1}^\infty \frac{1}{F_{9^n n} c_n} + \sum_{n=1}^\infty \frac{1}{F_{9^n n+1} c_n}$$

is a transcendental number. To see this, write

$$x \sum_{n=1}^\infty \frac{1}{F_{9^n n} c_n} + \sum_{n=1}^\infty \frac{1}{F_{9^n n+1} c_n} = \sum_{n=1}^\infty \frac{F_{9^n n+1} x + F_{9^n n}}{F_{9^n n} F_{9^n n+1} c_n}.$$

Aiming to use Theorem 1.4, pick $x_2 = x$, $\beta = 1/2$, $y = 1$, and any $0 < \alpha < 1 < \varepsilon$. Suppose $x \neq \bar{\varphi}$. The transcendence follows if we can find $\zeta \in \mathbb{C}$ such that $\Re(\zeta(F_{9^n n+1}x + F_{9^n n})) > 0$ for all sufficiently large $n$. If $\Im(x) \neq 0$, pick $\zeta = -i\Im(x)$. Otherwise, pick $\zeta = x - \bar{\varphi}$, as then

$$\zeta(F_{9^n n+1}x + F_{9^n n}) = F_{9^n n}(x - \bar{\varphi})\left(x + \frac{F_{9^n n - 1}}{F_{9^n n}}\right) \tag{13}$$

and $\lim_{n\to\infty} F_{9^n n-1}/F_{9^n n} = 1/\varphi = -\bar{\varphi}$ ensure that $\zeta(F_{9^n n+1}x + F_{9^n n})$ is a positive real number when $n$ is sufficiently large, and we are done.

This leaves us with the case of $x = \bar{\varphi}$, where we have

$$F_{9^n n+1}\bar{\varphi} + F_{9^n n} = \bar{\varphi}^{9^n n+1}.$$

While we have no hope of getting $\Re(\zeta\bar{\varphi}^{9^n n+1}) > 0$ for all large $n$, Remark 1.5 allows us to ignore this if we can show that

$$s_N := \sum_{n=1}^{N-1} \frac{\bar{\varphi}^{9^n n+1}}{F_{9^n n} F_{9^n n+1}}.$$

does not take the same value for infinitely many $N$. To see this, we let $M > N$ and use the converse triangle inequality to find

$$|s_N - s_M| = \left|\sum_{n=N}^{M} \frac{\bar{\varphi}^{9^n n+1}}{F_{9^n n} F_{9^n n+1}}\right| > 0,$$

for all sufficiently large $N$, and the example is complete.

The next example shows Theorem 1.4 is not easily replaced by Theorem 1.7 in the above example. For this purpose, we will need a simple lemma, which will be proven in Sect. 3, right after Lemma 3.6.

**Lemma 2.2** *Let $x$ be a fixed non-zero algebraic number, and let $a, b \in \mathbb{Z}$. Then there is a constant $C > 0$, depending only on $x$, so that $|a + bx| \geq C \max\{|a|, |b|, 1\}^{-2 \deg x}$ when $a + bx \neq 0$.*

*Example 2.3* In Example 2.1, Theorem 1.7 would not have been immediately applicable on the sequences that appears when $x \neq \bar{\varphi}$ is quadratic irrational. To see this, notice that we must have

$$a_n = F_{9^n n} F_{9^n n+1} \tilde{a}_n \quad \text{and} \quad b_n = (F_{9^n n+1} x + F_{9^n n}) \tilde{a}_n,$$

for some suitable sequence of $\tilde{a}_n \in \overline{\mathbb{Q}}$ such that $\kappa a_n$ and $\kappa b_n$ are all algebraic integers for some fixed $\kappa \in \mathbb{N}$. If $\tilde{a}_n \notin \mathbb{Q}(x)$ for some $n$, then we get $d \geq 4$ and so

$$\frac{d^2(y_1 + y_2)}{1 - \beta} + 1 \geq 17 > 9,$$

making Theorem 1.7 inapplicable, so we assume $\tilde{a}_n \in \mathbb{Q}(x)$ for all $n$. Due to Eq. (13), it follows that

$$y_2 \geq \beta \geq \frac{\limsup_{n \to \infty} \log |b_n|}{\limsup_{n \to \infty} \log |a_n|} = \frac{1 + c}{2 + c}, \quad \text{where } c = \limsup_{n \to \infty} \frac{\log |\tilde{a}_n|}{\log F_{9^n n}}.$$

Note that we need $c > -2$ in order to have a $y_1$ that satisfies the bound on $a_{i,n}$ for all $n$. By Remark 1.9, we may assume $D = 2$, $x_1 = 1$, and $x_2 = x$. Writing $\tilde{a}_n = \tilde{a}_{1,n} + \tilde{a}_{2,n} x$ with $\tilde{a}_{1,n}, \tilde{a}_{2,n} \in \mathbb{Q}$, we obtain from Lemma 2.2 that

$$|\tilde{a}_n| \geq C \max\{|\tilde{a}_{1,n}|, |\tilde{a}_{2,n}|\}^{-4},$$

where $C > 0$ is a constant that depends only on $x$, and so

$$y_1 \geq \limsup_{n \to \infty} \frac{\log(F_{9^n n} F_{9^n n+1} \max\{|\tilde{a}_{1,n}|, |\tilde{a}_{2,n}|\})}{\log(F_{9^n n} F_{9^n n+1} |\tilde{a}_n|)} \geq \frac{2 - \frac{c}{4}}{2 + c}.$$

If $c \leq -1$, then

$$\frac{d^2(y_1 + y_2)}{1 - \beta} + 1 \geq \frac{4\left(\frac{2 - c/4}{2 + c} + 0\right)}{1 - 0} + 1 = \frac{8 - c}{2 + c} + 1 \geq 9 + 1 > 9.$$

while assuming $c > -1$ yields

$$\frac{d^2(y_1 + y_2)}{1 - \beta} + 1 \geq \frac{4\left(\frac{2 - c/4}{2 + c} + \frac{1 + c}{2 + c}\right)}{1 - \frac{1 + c}{2 + c}} + 1 = 12 + 3c + 1 > 9.$$

This shows that Theorem 1.7 is not immediately applicable and concludes the example.

In the remaining examples, we fix $K = \mathbb{Q}(\sqrt{5}) = \mathbb{Q}(\varphi)$, $D = d = 2$, and $x_1 = 1$, while the value of $x_2$ will be chosen as either $\varphi$ or $\bar{\varphi}$. The aim of these examples is to show some more simple applications of the theorems while also highlighting the differences in applicability.

*Example 2.4* The sequence $\{n^{5^n}\varphi^n\}_{n=1}^\infty$ is transcendental. This follows from Theorem 1.6 with $x_2 = \varphi$, $\beta = 0$, $y = \eta_1 = \eta_2 = 1$, and arbitrary $\alpha, \delta, \varepsilon \in (0, 1)$. Alternatively, one could apply Theorem 1.4 with $x_2 = \bar{\varphi}$, $\beta = 0$, $y = 1$, and any $\alpha, \varepsilon \in (0, 1)$ upon rewriting to

$$n^{5^n}\varphi^n = \frac{n^{5^n}}{(-\bar{\varphi})^n} = \frac{n^{5^n}}{(-1)^n(F_n\bar{\varphi} + F_{n-1})}.$$

*Example 2.5* The sequence $\{\varphi^{7^n}\}_{n=1}^\infty$ is transcendental. This follows from Theorem 1.6 with $x_2 = \varphi$, $\beta = \eta_1 = 0$, and $\eta_2 = y = 1$. Note that if we wished to immediately apply Theorem 1.4, we would be left with

$$b_n = a_n\varphi^{-7^n} = a_n(F_{7^n-1} - F_{7^n}\varphi^{-1}).$$

Now, if $a_n \leq F_{7^n}$, we have $\beta = 0$ but must take $y \geq 2$, which means that the limsup criterion is not satisfied since $d^2y/(1-\beta) + 1 \geq 9 > 7$. If $a_n \geq F_{7^n}$ is large enough, we may achieve $y = 2 - \gamma$ for some $0 \leq \gamma < 1$, but then it is easily shown that we get $\beta \geq \gamma$, so that

$$\frac{d^2y}{1-\beta} + 1 \geq 9 > 7,$$

and the divergence criterion remains unsatisfied.

*Example 2.6* The sequence $\{F_{9^n n+1}\varphi^{9^n n}/F_{9^n n}\}_{n=1}^\infty$ is transcendental. This follows by rewriting into

$$\frac{F_{9^n n+1}\varphi^{9^n n}c_n}{F_{9^n n}} = \frac{F_{9^n n+1}}{F_{9^n n}(-\bar{\varphi})^{9^n n}} = \frac{F_{9^n n+1}}{F_{9^n n}(-1)^n(F_{9^n n-1} + F_{9^n n}\bar{\varphi})},$$

and then applying Theorem 1.4 with $x_2 = \bar{\varphi}$, $\beta = 0$, and $y = 2$ or, alternatively, Theorem 1.7 with $x_2 = \bar{\varphi}$, $\beta = 0$, $\eta_1 = \eta_2 = y_1 = 1$, and $y_2 = 2$. On the other hand, Theorem 1.6 is not immediately applicable, as can be seen through similar arguments to those in Example 2.5.

In our final example, below, we consider a sequence where only Theorem 1.7 is applicable. This will also serve as an example that we may encounter

$$|a_n|\left(\frac{\eta_2+d((d-1)y_1+y_2+\eta_2-\eta_1)+\delta}{1-\beta}+1\right)^{-n} > |a_n|\left(\frac{d^2(y_1+y_2)}{1-\beta}+1\right)^{-n}$$

for some sequences, though it should be mentioned that there appears to be no connection between this inequality and the applicability of the other theorems.

*Example 2.7* The sequence $\{\varphi^{2\cdot 14^n}/(F_{14^n} + \varphi)\}_{n=1}^\infty$ is transcendental. This follows from Theorem 1.7 by taking $\eta_1 = 0$, $\beta = y_1 = 1/2$, and $\eta_2 = y_2 = 1$. Here, neither one of Theorems 1.4 and 1.6 is immediately applicable, as seen through similar arguments to those in Example 2.5.

## 3 Preliminaries

A central tool for proving the main results is Schmidt's Subspace Theorem [8], below.

**Theorem 3.1** (Schmidt) *Let $L_1, \ldots, L_d$ be $\mathbb{Q}$-linearly independent linear forms in d variables with algebraic coefficients. For any $\delta > 0$, there exists a finite collection of proper subspaces $T_1, \ldots, T_w \subsetneq \mathbb{Q}^d$ such that any $x \in \mathbb{Z}^d$ with*

$$|L_1(x) \cdots L_n(x)| < |x|^{-\delta}$$

*is contained in $\bigcup_{i=1}^{w} T_i$.*

This theorem will be used together with the following lemma, which is found in a paper by Hančl, Nair, and Šustek [7]. The $M$ used in the present version of the lemma equals 1 plus the $M$ used in the original. [7] also had additional assumptions for the lemma (such as the $a_n$ being integers and $M$ having a greater lower bound), but those assumptions were never used in the proof and were only there for the sake of the main theorem of that paper.

**Lemma 3.2** *Consider real numbers $\varepsilon > 0$, $0 < \alpha < 1 \leq M$, and $\beta \in [0, \frac{1+\varepsilon}{\varepsilon})$. Let $\{a_n\}_{n=1}^{\infty}$ be a sequence of positive integers such that*

$$n^{1+\varepsilon} \leq a_n \leq a_{n+1} \quad \text{and} \quad \limsup_{n \to \infty} a_n^{\left(\frac{M}{1-\beta}+1\right)^{-n}} = \infty.$$

*Let $x_1, \ldots, x_d \in \mathbb{C}$, let $y_1, \ldots, y_d \geq 1$, and let $\{b_n\}_{n=1}^{\infty}$ be a sequence of complex numbers with $b_n = \sum_{i=1}^{d} x_i b_{i,n}$ for some $b_{1,n}, \ldots, b_{d,n} \in \mathbb{Z}$, such that for all sufficiently large n and each $i = 1, \ldots, d$,*

$$b_n \leq a_n^\beta 2^{\log_2^\alpha a_n}$$

*and*

$$|b_{i,n}| \leq |a_n|^{y_i} 2^{\log_2^\alpha |a_n|}.$$

*Finally, let $\{c_n\}_{n=1}^{\infty}$ be a sequence of positive integers. Then there is a positive real number $E > 0$ such that the inequality*

$$\left| \sum_{n=1}^{\infty} \frac{b_n}{a_n c_n} - \frac{\sum_{i=1}^{d} p_i x_i}{q} \right| < \frac{1}{(\log_2^2 q) 2^{d \log_2^{(1+2\alpha)/3} q} q^M} \tag{14}$$

*has infinitely many solutions $(p_1, \ldots, p_d, q) \in \mathbb{Z}^d \times \mathbb{N}$ satisfying*

$$|p_i| \leq E 2^{\log_2^{(1+2\alpha)/3} q} q^{y_i}, \quad \text{for all } i = 1, \ldots, d. \tag{15}$$

Another central tool for the proofs is the following lemma, which is a slight strengthening of another lemma from [7] and follows from a much similar proof. For clarity, we will go through the proof in Sect. 6.

**Lemma 3.3** *Let $\varepsilon > 0$, $0 < \alpha < 1 \leq M$, and $\beta \in [0, \frac{1+\varepsilon}{\varepsilon})$. Let $\{a_n\}_{n=1}^{\infty}$ and $\{b_n\}_{n=1}^{\infty}$ be sequences of positive real numbers such that for all sufficiently large n,*

$$n^{1+\varepsilon} \leq a_n \leq a_{n+1}, \qquad b_n \leq a_n^\beta 2^{\log_2^\alpha a_n},$$

*and*

$$\limsup_{n\to\infty} a_n^{\left(\frac{M}{1-\beta}+1\right)^{-n}} = \infty. \tag{16}$$

*Let $0 < c < 1$ be fixed. Then*

$$\liminf_{N\to\infty} \left(2^{N^2 \log_2^c a_{N-1}} \left(\prod_{n=1}^{N-1} a_n^M\right) \sum_{n=N}^{\infty} \frac{b_n}{a_n}\right) = 0.$$

We now present a few notions from algebraic number theory that will be relevant in the proofs of the main theorems. Let $a$ be an algebraic number with minimal polynomial $\sum_{i=0}^{d} c_i X^i$ over the integers ($c_d > 0$). The leading coefficient, $c_d$, is also called the *denominator* of $a$, since $c_d a$ is an algebraic integer while $c' a$ is not for any rational integer $0 < c' < c_d$. By rewriting the minimal polynomial of $a$ as $c_d \prod_{i=1}^{d}(X - a_i)$ instead, the Mahler measure of $a$ is defined as

$$M(a) := c_d \prod_{i=1}^{d} \max\{1, |a_i|\}.$$

Surprisingly closely related to this is the Weil height, which we define as

$$H(\alpha) := \prod_{v \in M_K} \max\{1, |a|_v\}^{[K_v:\mathbb{Q}_v]/[K:\mathbb{Q}]},$$

where $K$ is any number field containing $a$, $M_K$ denotes the set of places of $K$, $K_v$ is the local field of $K$ at $v$, and $[K : L]$ denotes the degree of a field extension $K \supseteq L$. This does not depend on the choice of $K$ (see [9] for a proof). We will compare and estimate the house, Mahler measure, and Weil height using the following classical results.

**Lemma 3.4** *Let $a$ be an algebraic number with denominator $c_d$. Then*

$$H(a)^d = M(a) \leq |c_d| \max\{\overline{|a|}^d, 1\}.$$

*Proof* The inequality is a trivial consequence of the definitions. For the equality, see [9, Lemma 3.10]. □

**Lemma 3.5** *Let $a, b \in \overline{\mathbb{Q}}$ with $a \neq 0$. Then*

$$H(a+b) \leq 2H(a)H(b), \quad H(ab) \leq H(a)H(b),$$
$$H(1/a) = H(a).$$

*Proof* See [9]. □

**Lemma 3.6** (Liouville Inequality) *Let $a, b$ be non-conjugate algebraic numbers. Then*

$$|a - b| \geq \left(2H(a)H(b)\right)^{-\deg(a)\deg(b)}.$$

*Proof* This can be extracted from [10, Theorem A.1]. □

*Proof of Lemma* 2.2 Pick $C = \min\{(2H(x))^{-\deg x}, |x|\}$. If $a$ or $b$ is 0, the statement is trivially true. If $ab \neq 0$, then we may apply Lemmas 3.6 and 3.5 to conclude

$$|a + bx| \geq (2H(a)H(-bx))^{-\deg(a)\deg(-bx)} \geq (2H(a)H(b)H(x))^{-\deg(x)}$$
$$= (2H(x))^{-\deg(x)}|ab|^{-\deg x} \geq C \max\{|a|, |b|\}^{-\deg x}.$$

$\square$

In order to prove Theorem 1.7, we will need a different version of Lemma 3.2, the proof of which will use some elementary of Galois theory. Recall that a field extension $K \supseteq \mathbb{Q}$ is called a Galois extension if all irreducible polynomials over $\mathbb{Q}$ that have a root in $K$ split into linear factors over $K$. Note that this is equivalent to $K$ being closed under conjugation. We use $\mathrm{Gal}(K)$ to denote the associated Galois group, i.e., the field automorphisms on $K$ that preserve $\mathbb{Q}$. Recall that for any finite field extension $L \supseteq \mathbb{Q}$, $L$ has a unique finite field extension $K \supseteq L$ of minimal degree such that $K \supseteq \mathbb{Q}$ is Galois (see, e.g., Theorem 11.6 of [11]). This also implies the below lemma, which will be relevant for the proofs of both Theorems 1.4 and 1.7.

**Lemma 3.7** *Let $a_1, \ldots, a_d \in \overline{\mathbb{Q}}$. Then there is a constant C, depending only on $a_1, \ldots, a_d$, so that for any $(c_1, \ldots, c_d) \in \mathbb{Q}^d$,*

$$\overline{|c_1 a_1 + c_2 a_2 + \cdots + c_d a_d|} \leq C \max_{1 \leq i \leq d} |c_i|.$$

*Proof* Let $K \supseteq \mathbb{Q}$ be the smallest Galois extension of $\mathbb{Q}$ containing $a_1, \ldots, a_d$. Since conjugation is a field automorphism on $\overline{\mathbb{Q}}$, and $K$ is closed under conjugation, we find

$$\overline{\left|\sum_{i=1}^d c_i a_i\right|} \leq \max_{\psi \in \mathrm{Gal}(K)} \left|\psi\left(\sum_{i=1}^d c_i a_i\right)\right| = \max_{\psi \in \mathrm{Gal}(K)} \left|\sum_{i=1}^d c_i \psi(a_i)\right|.$$

The proof is then completed by an application of the triangle inequality,

$$\overline{\left|\sum_{i=1}^d c_i a_i\right|} \leq \left(d \max_{\psi \in \mathrm{Gal}(K)} \max_{0 \leq i \leq d} |\psi(a_i)|\right) \max_{0 \leq i \leq d} |c_i|.$$

$\square$

## 4 Proof of Theorem 1.4

We will first prove the below result, which is inspired by the main theorem of [7].

**Theorem 4.1** *Let $d \in \mathbb{N}$ be a positive integer, and consider real numbers $\delta, \varepsilon > 0$, $0 < \alpha < 1 \leq y_1, \ldots, y_n$, and $\beta \in [0, \frac{\varepsilon}{1+\varepsilon})$. Let furthermore $x_1, \ldots, x_d$ be algebraic numbers, and let $\{a_n\}_{n=1}^\infty$ be a sequence of positive integers that satisfy inequality (2) and*

$$\limsup_{n \to \infty} a_n^{\left(\frac{1+\sum_{i=1}^d y_i + \delta}{1-\beta} + 1\right)^{-n}} = \infty.$$

*Let $\{b_n\}_{n=1}^\infty$ be a sequence of non-zero numbers given by $b_n = \sum_{i=1}^d b_{i,n} x_i$ where $b_{i,n} \in \mathbb{Z}$ and such that the inequalities of inequality (3) are satisfied for n sufficiently large and each $i = 1, \ldots, d$. Let $\{c_n\}_{n=1}^\infty$ be a sequence of positive integers. Then the number $\sum_{n=1}^\infty \frac{b_n}{a_n c_n}$ is either transcendental or a $\mathbb{Q}$-linear combination of $x_1, \ldots, x_d$.*

The main difference between the proof of this theorem and that of the corresponding one in [7] lies in the below application of Theorem 3.1, which replaces [7, Lemma 7].

**Lemma 4.2** *Let $x_1, \ldots, x_d, s$ be algebraic numbers such that $s$ is $\mathbb{Q}$-linearly independent of $x_1, \ldots, x_d$, and let $\delta > 0$. Then the inequality,*

$$\left| qs - \sum_{i=1}^{d} p_i x_i \right| \prod_{i=1}^{d} \max\{1, |p_i|\} < q^{-\delta}, \tag{17}$$

*has only finitely many solutions $(p_1, \ldots, p_d, q) \in \mathbb{Z}^d \times \mathbb{N}$.*

*Proof* This will be proven by induction, using the convention that linear independence of the empty set is equivalent to being non-zero. Let $S$ denote the set of solutions $(p_1, \ldots, p_d, q) \in (\mathbb{Z} \setminus \{0\})^d \times \mathbb{N}$ to inequality (17). For $d = 0$ or $S = \emptyset$, $S$ is clearly finite, so suppose $d > 0$, $S \neq \emptyset$, and that the lemma is true for $d' = d - 1$.

Note that all elements of $S$ satisfy

$$\left| qs - \sum_{i=1}^{d} p_i x_i \right| \prod_{i=1}^{d} |p_i| < q^{-\delta}.$$

By Theorem 3.1, there is a finite collection of proper subspaces $T_1, \ldots, T_w \subsetneq \mathbb{Q}^{d+1}$ such that $S \subseteq \bigcup_{l=1}^{w} T_l$. Write $S_l = S \cap T_l$ for $l = 1, \ldots, w$, and let $1 \leq l \leq w$ such that $S_l \neq \emptyset$. Then $T_l$ contains an element with a non-zero $q$-entry. Since $\dim T_l \leq d$, it follows that there is a $j$ such that the $p_j$-entry is given as a fixed linear combination of the remaining entries for all elements in $T_l$. By renumbering if necessary, we may assume that $j = d$ and then pick $r_1, \ldots, r_d \in \mathbb{Q}$ such that $p_d = r_1 p_1 + \cdots + r_{d-1} p_{d-1} + r_d q$ for all $(p_1, \ldots, p_d, q) \in T_l$. For elements of $S_l$, inequality (17) now reduces to

$$q^{-\delta} > \left| q(s - x_d r_d) - \sum_{i=1}^{d-1} p_i (x_i + x_d r_i) \right| \prod_{i=1}^{d} \max\{1, |p_i|\}$$
$$\geq \left| q(s - x_1 r_1) - \sum_{i=1}^{d-1} p_i (x_i + x_1 r_i) \right| \prod_{i=1}^{d-1} \max\{1, |p_i|\}.$$

Hence, $S_l$ is finite by induction. Since $S = S_1 \cup \cdots \cup S_w$, this completes the proof. □

*Proof of Theorem 4.1* Put $M = 1 + \sum_{i=1}^{d} y_i + \delta$. By Lemma 3.2, there are infinitely many $(p_1, \ldots, p_d, q) \in \mathbb{Z}^d \times \mathbb{N}$ satisfying both inequalities (14) and (15), where we may take $E$ to be rational. Rewriting inequality (14) using the above choice of $M$, we find

$$\left| q \sum_{n=1}^{\infty} \frac{b_n}{a_n c_n} - \sum_{i=1}^{d} p_i x_i \right| \prod_{i=1}^{d} \left( q^{y_i} 2^{\log_2^{(1+2\alpha)/3} q} \right) < q^{-\delta},$$

and it follows from inequality (15) that

$$\left| qE^{-d} \sum_{n=1}^{\infty} \frac{b_n}{a_n c_n} - \sum_{i=1}^{d} p_i E^{-d} x_i \right| \prod_{i=1}^{d} \max\{1, |p_i|\} < q^{-\delta}.$$

Lemma 4.2 then implies that $E^{-d} \sum_{n=1}^{\infty} \frac{b_n}{a_n c_n}$ cannot both be algebraic and $\mathbb{Q}$-linearly independent of $E^{-d} x_1, \ldots, E^{-d} x_k$. Since $E$ is rational, this completes the proof. □

To finish the proof of Theorem 1.4, we will use the below proposition to ensure that $\sum_{n=1}^{\infty} \frac{b_n}{a_n c_n}$ is indeed $\mathbb{Q}$-linearly independent of $x_1, \ldots, x_d$.

**Proposition 4.3** *Let $d, \tilde{d} \in \mathbb{N}$, let $K$ be a number field of degree $d$, let $x_1, \ldots, x_D \in K$, and consider real numbers $\varepsilon > 0$, $0 < \alpha < 1$, $\beta \in [0, \frac{\varepsilon}{1+\varepsilon})$, and $y \geq 1$. Let $\{a_n\}_{n=1}^{\infty}$ be a sequence of positive integers that satisfy inequality (2) and*

$$\limsup_{n \to \infty} a_n^{\left(\frac{d\tilde{d}y}{1-\beta}+1\right)^{-n}} = \infty.$$

*Let $\{b_n\}_{n=1}^{\infty}$ be a sequence of non-zero numbers given by $b_n = \sum_{i=1}^{d} b_{i,n} x_i$ where $b_{i,n} \in \mathbb{Z}$ and such that inequalities (3) and (4) are satisfied for all sufficiently large $n$. Let $\{c_n\}_{n=1}^{\infty}$ be a sequence of positive integers. Then the number $\sum_{n=1}^{\infty} \frac{b_n}{a_n c_n}$ has degree strictly greater than $\tilde{d}$.*

*Proof* Let $\sigma : \mathbb{N} \to \mathbb{N}$ be a bijection such that $A_n = a_{\sigma(n)} c_{\sigma(n)}$ is increasing, and put $B_n = b_{\sigma(n)}$ and $B_{i,n} = B_{i,\sigma(n)}$. Since clearly $A_n \geq a_n$ for all $n$, we get that $\{A_n\}_{n=1}^{\infty}$, $\{B_{i,n}\}_{n=1}^{\infty}$, and $\{B_n\}_{n=1}^{\infty}$ satisfy all assumptions of the proposition. For the remainder of the proof, we may therefore assume that $c_n = 1$ for all $n$.

Assume towards contradiction that $\deg\left(\sum_{n=1}^{\infty} \frac{b_n}{a_n}\right) \leq \tilde{d}$, and write

$$s = \sum_{n=1}^{\infty} \frac{b_n}{a_n}, \qquad s_N = \sum_{n=1}^{N-1} \frac{b_n}{a_n}.$$

Note that $\deg(s_N) \leq d$, and let $c$ denote the least common multiple of the denominators of $x_1, \ldots, x_D$. Then the denominator of $s_N$ is at most $c \prod_{n=1}^{N-1} a_n$. By Lemmas 3.4 and 3.7,

$$H(s_N) \leq c \left(\prod_{n=1}^{N-1} a_n\right) \overline{|s_N|} \leq C_1 \left(\prod_{n=1}^{N-1} a_n\right) \max_{1 \leq i \leq D} \left|\sum_{n=1}^{N-1} \frac{b_{i,n}}{a_n}\right|, \tag{18}$$

where $C_1 > 0$ is some sufficiently large constant depending only on $x_1, \ldots, x_D$. Since $y \geq 1$, the triangle inequality and inequality (3) then imply

$$\left|\sum_{n=1}^{N-1} \frac{b_{i,n}}{a_n}\right| \leq \sum_{n=1}^{N-1} a_n^{y-1} 2^{\log_2^{\alpha} a_n} \leq N 2^{\log_2^{\alpha} a_N} a_{N-1}^{y-1},$$

for all sufficiently large $N$. Applying this to inequality (18) and once again using that $y \geq 1$, we obtain that when $N$ is sufficiently large,

$$H(s_N) \leq C_1 \left(\prod_{n=1}^{N-1} a_n\right) N a_{N-1}^{y-1} 2^{\log_2^{\alpha} a_{N-1}} \leq 2^{(3d\tilde{d})^{-1} N^2 \log_2^{\alpha} a_{N-1}} \prod_{n=1}^{N-1} a_n^y. \tag{19}$$

When $N$ grows large, inequality (4) makes $\Re(\zeta s_N)$ strictly increasing. Since $s$ has only finitely many conjugates, $s$ and $s_N$ can thus only be conjugate numbers for finitely many $N$. When $N$ is sufficiently large, it therefore follows from the triangle inequality and Lemma 3.6 that

$$\sum_{n=N}^{\infty} \left|\frac{b_n}{a_n}\right| \geq |s - s_N| \geq (2H(s)H(s_N))^{-\deg(s)\deg(s_n)},$$

and so, recalling that $\deg s \leq \tilde{d}$ and $\deg s_n \leq d$ while applying inequality (19),

$$\sum_{n=N}^{\infty} \left|\frac{b_n}{a_n}\right| \geq \left(2H(s) 2^{(3d\tilde{d})^{-1} N^2 \log_2^{\alpha} a_{N-1}} \prod_{n=1}^{N-1} a_n^y\right)^{-d\tilde{d}} > \sqrt{2}^{N^2 \log_2^{\alpha} a_{N-1}} \prod_{n=1}^{N-1} a_n^{-d\tilde{d}y}$$

We conclude that for all large enough $N$,

$$2^{N^2 \log_2^\alpha a_{N-1}} \left( \prod_{i=1}^{N-1} a_n^{d\tilde{d}y} \right) \sum_{n=N}^\infty \left| \frac{b_n}{a_n} \right| > \sqrt{2}^{N^2 \log_2^\alpha a_{N-1}},$$

which contradicts Lemma 3.3 and thus completes the proof. □

*Proof of Theorem 1.4* By Remark 1.9, we may assume that $D = d$ and that $x_1, \ldots, x_d$ forms a $\mathbb{Q}$-linear basis of $K$. The irrationality statement is then simply Proposition 4.3 with $\tilde{d} = 1$, while the transcendence statement follows from Theorem 4.1 and Proposition 4.3 with $\tilde{d} = d$. □

## 5 Proof of Theorem 1.7

Our first step in proving Theorem 1.7 will be to prove the below parallel result to Theorem 4.1.

**Theorem 5.1** *Let $K$ be a number field with $\mathbb{Q}$-linear basis $x_1, \ldots, x_d$, and consider real numbers $\alpha, \varepsilon > 0$, $\beta, \eta_1, y_1 \geq 0$, and $\eta_2, y_2 \geq 1$ such that $\alpha < 1$, $\beta < \varepsilon/(1 + \varepsilon)$, and $\eta_1 \leq (d-1)y_1 + y_2$. Let $\{a_n\}_{n=1}^\infty$ and $\{b_n\}_{n=1}^\infty$ be non-zero sequences in $K$ given by $a_n = \sum_{i=1}^d a_{i,n} x_i$ and $b_n = \sum_{i=1}^d b_{i,n} x_i$ where $a_{i,n}, b_{i,n} \in \mathbb{Z}$ such that*

$$\limsup_{n \to \infty} |a_n|^{\left(\frac{\eta_2 + d((d-1)y_1 + y_2 - \eta_1 + \eta_2) + \delta}{1-\beta} + 1\right)^{-n}} = \infty$$

*and inequalities (6), (7), (8), (9), and (10) are satisfied for each $i = 1, \ldots, d$ and all sufficiently large $n$. Then the number $\sum_{n=1}^\infty \frac{b_n}{a_n c_n}$ is either transcendental or a $\mathbb{Q}$-linear combination of $x_1, \ldots, x_d$.*

This is not quite as neat as Theorem 4.1, and the reason for this is to be found in Lemma 3.2. As part of its proof in [7], the authors write $b_n/a_n$ as a $\mathbb{Q}$-linear combination of $x_1, \ldots, x_d$, which is fairly elegantly done when $a_n$ is rational and less so when $a_n$ may be irrational. Using the Galois theory mentioned by the end of Sect. 3 to make the corresponding modifications, we reach the below lemma.

**Lemma 5.2** *Using the notation and assumptions of Theorem 5.1, the inequality*

$$\left| \sum_{n=1}^\infty \frac{b_n}{a_n c_n} - \frac{\sum_{i=1}^d p_i x_i}{q} \right| < \frac{1}{2^{d \log_2^{(1+2\alpha)/3} q} q^{d+1+d\frac{(d-1)y_1+y_2-\eta_1}{\eta_2}}} \quad (20)$$

*has infinitely many solutions $(p_1, \ldots, p_d, q) \in \mathbb{Z}^d \times \mathbb{N}$ satisfying*

$$|p_i| \leq 2^{\log_2^{(1+2\alpha)/3} q} q^{1 + \frac{(d-1)y_1+y_2-\eta_1}{\eta_2}}, \quad \text{for all } i = 1, \ldots, d. \quad (21)$$

*Proof* Let $\sigma : \mathbb{N} \to \mathbb{N}$ be a bijection such that the sequence $\{A_n\}_{n=1}^\infty$, given by $A_n = a_{\sigma(n)} c_{\sigma(n)}$, is of increasing modulus. Put $B_{i,n} = b_{i,\sigma(n)}$, $B_n = b_{\sigma(n)}$, and $R_n = r_{\sigma(n)}$. Since these new sequences satisfy the hypothesis of the lemma, we may assume without loss of generality that $c_n = 1$ for all $n$.

Let $\tilde{K} \supseteq \mathbb{Q}$ be the smallest Galois extension of $\mathbb{Q}$ with $K \subseteq \tilde{K}$. Pick $x_{d+1}, \ldots, x_D$ such that $x_1, \ldots, x_D$ is a $\mathbb{Q}$-linear basis of $\tilde{K}$, and let $\pi_i$ denote the $i$'th coordinate map in this basis. Note that $\pi_i(a_n) = a_{i,n}$ and $\pi_i(b_n) = b_{i,n}$ when $1 \leq i \leq d$. Pick $c > 0$ and $\kappa \in \mathbb{N}$ so

that

$$c \geq \left| \pi_i \left( \prod_{k=1}^{d'} g_k(x_{j_k}) \right) \right| \quad \text{and} \quad \kappa \pi_i \left( \prod_{k=1}^{d'} g_k(x_{j_k}) \right) \in \mathbb{Z},$$

for all $d' \mid d$, all $i, j_1, \ldots, j_{d'} \in \{1, \ldots, d\}$, and all $g_1, \ldots, g_{d'} \in \mathrm{Gal}(\tilde{K})$.

For each $n \in \mathbb{N}$, pick $\tilde{a}_n \in \mathbb{Z}$ of minimal modulus such that

$$\tilde{a}_n r_n \mathcal{N}(a_n/r_n) \geq |a_n|^{\eta_2},$$

and note by inequality (8) that then

$$\tilde{a}_n r_n \mathcal{N}(a_n/r_n) \leq 2|a_n|^{\eta_2} 2^{\log_2^\alpha |a_n|}. \tag{22}$$

Write $d_n = \deg a_n$ and pick $g_1, \ldots, g_{d_n} \in \mathrm{Gal}(\tilde{K})$ such that $(g_k(a_n))_{k=1}^{d_n}$ runs through all $d_n$ conjugates of $a_n$, with $g_1(a_n) = a_n$. It follows that

$$\pi_i \left( \kappa r_n \mathcal{N}\left(\frac{a_n}{r_n}\right) \frac{b_n}{a_n} \right) = \kappa \pi_i \left( \left( \sum_{j=1}^d b_{j,n} x_j \right) \prod_{k=2}^{d_n} \sum_{j=1}^d \frac{a_{j,n}}{r_n} g_k(x_j) \right)$$

$$= \sum_{j_1=1}^d b_{j_1,n} \sum_{j_2=1}^d \frac{a_{j_2,n}}{r_n} \cdots \sum_{j_{d_n}=1}^d \frac{a_{j_{d_n},n}}{r_n} \kappa \pi_i \left( \prod_{k=1}^{d_n} g_k(x_{j_k}) \right)$$

must be an integer by choice of $\kappa$ since $d_n \mid d$. Define

$$q_N := \kappa \prod_{n=1}^{N-1} \tilde{a}_n r_n \mathcal{N}(a_n/r_n) \quad \text{and} \quad p_{i,N} := \pi_i \left( q_N \sum_{n=1}^{N-1} \frac{b_n}{a_n} \right),$$

and note that $q_N \in \mathbb{N}$ and $p_{i,N} \in \mathbb{Z}$ by the above considerations. Set $M = \eta_2 + d((d-1)y_1 + y_2 - \eta_1 + \eta_2) + \delta$. The choice of $p_{i,N}$, the triangle inequality, and Lemma 3.3 show that for infinitely many $N$,

$$\left| \sum_{n=1}^\infty \frac{b_n}{a_n} - \frac{\sum_{i=1}^d p_{i,N} x_i}{q_N} \right| = \left| \sum_{n=1}^\infty \frac{b_n}{a_n} - \sum_{n=1}^{N-1} \frac{b_n}{a_n} \right| \leq \sum_{n=N}^\infty \left| \frac{b_n}{a_n} \right|$$

$$< 2^{-\log_2^{(1+\alpha)/2} \prod_{n=1}^{N-1} |a_n|} \prod_{n=1}^{N-1} |a_n|^{-M}.$$

Inequality (22), with the choices of $q_N$ and $M$, then implies

$$\left| \sum_{n=1}^\infty \frac{b_n}{a_n} - \frac{\sum_{i=1}^d p_{i,N} x_i}{q_N} \right| < 2^{-\log_2^{(1+\alpha)/2} \prod_{n=1}^{N-1} |a_n|} \prod_{n=1}^{N-1} (\tilde{a}_n r_n \mathcal{N}(a_n/r_n))^{-M/\eta_2}$$

$$\leq 2^{-\log_2^{(1+\alpha)/2} q_N} q_N^{-M/\eta_2}$$

$$= 2^{-\log_2^{(1+\alpha)/2} q_N} q_N^{-1-d-d\frac{(d-1)y_1+y_2-\eta_1}{\eta_2}}.$$

Since $\log_2^{(1+\alpha)/2} q_N \geq d \log_2^{(1+2\alpha)/3} q_N$ when $N$ (and thereby $q_N$) is sufficiently large, we conclude that inequality (20) is satisfied for $q = q_N$ and $p_i = p_{i,N}$, for infinitely many choices of $N$.

We now just need to check that inequality (21) is also satisfied for $q = q_N$ and $p_i = p_{i,N}$. We start by noting

$$\frac{b_n}{a_n} = \frac{1}{\mathcal{N}(a_n)^{d/d_n}} \left(\sum_{j=1}^d b_{j,n} x_j\right) \left(\sum_{j=1}^d a_{j,n} x_j\right)^{d/d_n - 1} \prod_{k=2}^{d_n} \left(\sum_{j=1}^d a_{j,n} g_k(x_j)\right)^{d/d_n}$$

$$= \frac{1}{\mathcal{N}_K(a_n)} \sum_{j_1=1}^d \cdots \sum_{j_d=1}^d b_{j_1,n} a_{j_2,n} \cdots a_{j_d,n} \prod_{k=1}^{d_n} \prod_{l=1}^{d/d_n} g_k(x_{j_{kd_n+l}}).$$

It then follows from the triangle inequality and the choice of $c$ that

$$\left|\pi_i\left(\frac{b_n}{a_n}\right)\right| \leq d^d c \frac{1}{\mathcal{N}(a_n)} \max_{1 \leq j \leq d} |b_{j,n}| \max_{1 \leq j \leq d} |a_{j,n}|^{d-1}.$$

Hence, by inequalities (7) and (9),

$$\left|\pi_i\left(\frac{b_n}{a_n}\right)\right| \leq d^d c \, |a_n|^{-\eta_1 + y_2 + y_1(d-1)} 2^{(d+1) \log_2^\alpha |a_n|}.$$

Recalling the choice of $p_{i,N}$, we now find

$$\left|\frac{p_{i,N}}{q_N}\right| \leq d^d c \sum_{n=1}^{N-1} |a_n|^{y_1(d-1) + y_2 - \eta_1} 2^{(d+1) \log_2^\alpha |a_n|}$$

$$\leq N |a_{N-1}|^{y_1(d-1) + y_2 - \eta_1} 2^{(d+2) \log_2^\alpha |a_N|}.$$

For all sufficiently large $N$, the choice of $q_N$ and inequality (6) lead to

$$\log_2^\alpha q_N \geq \log_2^\alpha \left((N-2)! |a_{N-1}|\right) \geq \frac{1}{d+3} \max\{\log_2 N, \log_2^\alpha |a_{N-1}|\}.$$

Since $\eta_1 \leq y_1(d-1) + y_2$, this means that

$$\left|\frac{p_{i,N}}{q_N}\right| \leq |a_{N-1}|^{y_1(d-1) + y_2 - \eta_1} 2^{\log_2^\alpha |q_N|} \leq |q_N|^{\frac{y_1(d-1) + y_2 - \eta_1}{\eta_2}} 2^{\log_2^\alpha |q_N|},$$

by choice of $\tilde{a}_N$ and $q_N$, and the proof is complete.    □

*Proof of Theorem 5.1* This follows in full parallel to the proof of Theorem 4.1, with $M = d(y_1(d-1) + y_2 + \eta_2 - \eta_1) + \delta$ and using Lemma 5.2 in place of Lemma 3.2. In the application of Lemma 4.2, replace $\delta$ with $\delta/\eta_2$.    □

We now just need a variant of Proposition 4.3 where we allow the $a_n$ to be irrational.

**Proposition 5.3** *Let $d, \tilde{d} \in \mathbb{N}$, let $K$ be a number field of degree $d$, let $x_1, \ldots, x_D \in K$, and consider real numbers $\varepsilon > 0$, $y_1 \geq 1 > \alpha > 0$, $\beta \in [0, \frac{\varepsilon}{1+\varepsilon})$, and $y_2 \geq \beta$. Let $\{a_n\}_{n=1}^\infty$ and $\{b_n\}_{n=1}^\infty$ be sequences of non-zero numbers given by $a_n = \sum_{i=1}^D a_{i,n} x_i$ and $b_n = \sum_{i=1}^D b_{i,n} x_i$ where $a_{i,n}, b_{i,n} \in \mathbb{Z}$ and satisfy inequality (6),*

$$n^{1+\varepsilon} \leq |a_n| \leq |a_{n+1}|,$$

*and*

$$\limsup_{n\to\infty} |a_n|^{\left(\frac{d\tilde{d}(y_1+y_2)}{1-\beta}+1\right)^{-n}} = \infty.$$

*Suppose for all sufficiently large n, each $i = 1\ldots,D$, and some fixed $\zeta \in \mathbb{C}$ that inequalities (9), (10), and (11) are satisfied. Let $\{c_n\}_{n=1}^{\infty}$ be a sequence of positive integers. Then the number $\sum_{n=1}^{\infty} \frac{b_n}{a_n c_n}$ has algebraic degree strictly greater than $\tilde{d}$.*

The main change from the proof of Proposition 4.3 is that we now use Lemma 3.5 a few times before using Lemmas 3.4 and 3.7 in the estimate of $H(s_N)$. This makes the proof closer to that in [6].

*Proof of Proposition 5.3* The proof is essentially the same as that of Proposition 4.3, except that the calculation starting with inequality (18) and ending with inequality (19) is replaced by

$$\begin{aligned}
H(s_N) &\leq 2^{N-2} \prod_{i=1}^{N-1} H(a_n) H(b_n) \leq 2^{N-2} \prod_{i=1}^{N-1} c^{2/d} \overline{|a_n|} \ \overline{|b_n|} \\
&\leq C^N \prod_{i=1}^{N-1} \max_{1\leq i \leq D} |a_{i,n}| \max_{1\leq i \leq D} |b_{i,n}| \leq C^N \prod_{i=1}^{N-1} |a_n|^{y_1+y_2} 2^{2\log_2^{\alpha}|a_n|}, \\
&\leq 2^{(3d\tilde{d})^{-1} N^2 \log_2^{\alpha} |a_{N-1}|} \prod_{i=1}^{N-1} |a_n|^{y_1+y_2},
\end{aligned}$$

using Lemmas 3.5, 3.4, 3.7 and the inequalities of inequality (9), where $C > 0$ is some sufficiently large constant that depends only on $x_1, \ldots, x_D$. The rest of the proof follows the exact same arguments as those for the proof of Proposition 4.3 but with $y$ replaced with $y_1 + y_2$.

□

*Proof of Theorem 1.7* By Remark 1.9, we may assume that $D = d$ and that $x_1, \ldots, x_d$ form a $\mathbb{Q}$-linear basis of $K$. Then the irrationality statement is identical to Proposition 5.3 with $\tilde{d} = 1$, while the transcendence statement follows from Theorem 5.1 and Proposition 5.3 with $\tilde{d} = d$.

□

## 6 Proof of Lemma 3.3

The proof of Lemma 3.3 closely follows the proof of [7, Lemma 5]. Similarly to Lemma 3.2, the below lemmas are taken from [7], in which the first three of them appear with additional assumptions that are never used in their proofs.

**Lemma 6.1** *$\{a_n\}_{n=1}^{\infty}$ and $\{b_n\}_{n=1}^{\infty}$ satisfy the assumptions of Lemma 3.3. Then there is a fixed number $0 < \gamma < 1$ such that for all sufficiently large N,*

$$\sum_{n=N}^{\infty} \frac{b_n}{a_n} \leq \frac{1}{a_N^{\gamma}}.$$

**Lemma 6.2** *Let $\beta$, $\{a_n\}_{n=1}^{\infty}$, and $\{b_n\}_{n=1}^{\infty}$ satisfy the assumptions of Lemma 3.3. Suppose that $a_n \geq 2^n$ for all sufficiently large n. Then there is a fixed number $0 < \Gamma < 1$ such that for all sufficiently large N,*

$$\sum_{n=N}^{\infty} \frac{b_n}{a_n} \leq \frac{2^{\log_2^{\Gamma} a_N}}{a_N^{1-\beta}}.$$

**Lemma 6.3** *Let $\beta$, $\{a_n\}_{n=1}^{\infty}$, and $\{b_n\}_{n=1}^{\infty}$ satisfy the assumptions of Lemma 3.3. Then there is a fixed number $0 < \Gamma < 1$ so that if $N \leq Q$ are sufficiently large and $a_n \geq 2^n$ for $n = N, \ldots, Q$, then*

$$\sum_{n=N}^{Q} \frac{b_n}{a_n} \leq \frac{2^{\log_2^{\Gamma} a_N}}{a_N^{1-\beta}}.$$

**Lemma 6.4** *Let $\{y_n\}_{n=1}^{\infty}$ be an unbounded sequence of positive real numbers. Then there are infinitely many positive integers N such that*

$$y_N > \left(1 + \frac{1}{N^2}\right) \max_{1 \leq n < N} y_n.$$

By a simple induction argument, we notice for $k < N$ and $\delta \geq 0$ that

$$(M+1+\delta)^N = (M+1+\delta)^{N-1} + (M+\delta)(M+1+\delta)^{N-1} = \cdots$$
$$= (M+1+\delta)^k + (M+\delta)\sum_{n=k}^{N-1}(M+1+\delta)^n. \quad (23)$$

This will be used both in the proof of Lemma 3.3 and for proving the below lemma, which is to be used together with Lemma 6.4 and Eq. (16).

**Lemma 6.5** *Let $\{a_n\}_{n=1}^{\infty}$ be a sequence of positive real numbers, and let k be a positive integer. Then for all $N > k$,*

$$\left(\max_{k \leq n < N} a_n^{\left(\frac{M}{1-\beta}+1\right)^{-n}}\right)^{\left(\frac{M}{1-\beta}+1\right)^N} > \prod_{n=k}^{N-1} a_n^{\frac{M}{1-\beta}}.$$

*Proof* We use equation (23) with $\delta = 0$ to find

$$\left(\max_{k \leq n < N} a_n^{\left(\frac{M}{1-\beta}+1\right)^{-n}}\right)^{\left(\frac{M}{1-\beta}+1\right)^N} > \left(\max_{k \leq n < N} a_n^{\left(\frac{M}{1-\beta}+1\right)^{-n}}\right)^{\frac{M}{1-\beta}\sum_{n=k}^{N-1}\left(\frac{M}{1-\beta}+1\right)^n}$$
$$\geq \prod_{n=k}^{N-1} \left(a_n^{\left(\frac{M}{1-\beta}+1\right)^{-n}}\right)^{\frac{M}{1-\beta}\left(\frac{M}{1-\beta}+1\right)^n} = \prod_{n=k}^{N-1} a_n^{\frac{M}{1-\beta}}.$$

□

*Proof of Lemma 3.3* We will split into three cases, depending on whether

$$\limsup_{n \to \infty} a_n^{\left(\frac{M}{1-\beta}+1+\delta\right)^{-n}} = \infty \quad (24)$$

holds for some fixed $\delta > 0$ and whether $a_n < 2^n$ infinitely often. To shorten notation, write

$$Z_N = 2^{N^2 \log_2^c a_{N-1}} \left( \prod_{n=1}^{N-1} a_n^M \right) \sum_{n=N}^{k_2} \frac{b_n}{a_n}.$$

**Case 1 (inequality (24) is satisfied for some $\delta > 0$)** Pick $0 < \gamma < 1$ as in Lemma 6.1, and let $z > 2$ be some sufficiently large number. Pick $k_1, k_2, N \in \mathbb{N}$ as follows. Let $k_2$ be the smallest integer such that

$$a_{k_2}^{\left(\frac{M}{1-\beta}+1+\delta\right)^{-k_2}} > z^{1/\gamma}, \tag{25}$$

let $k_1$ be the largest integer such that $k_1 < k_2$ and

$$a_{k_1} \leq z^{k_1}, \tag{26}$$

and let $N$ be the smallest number such that $N > k_1$ and

$$a_N^{\left(\frac{M}{1-\beta}+1+\delta\right)^{-N}} \geq z. \tag{27}$$

Note that $k_2 \geq N > k_1$ and that $k_1 \to \infty$ as $z \to \infty$. From the above choices of $k_1$ and $N$, it follows that $a_n < z^{(\frac{M}{1-\beta}+1+\delta)^n}$ when $k_1 \leq n < N$. Hence, by also applying Eq. (23)

$$\prod_{n=k_1}^{N-1} a_n < z^{\sum_{n=k_1}^{N-1} (\frac{M}{1-\beta}+1+\delta)^n} < z^{(\frac{M}{1-\beta}+\delta)^{-1} (\frac{M}{1-\beta}+1+\delta)^N},$$

while inequality (26) implies

$$\prod_{n=1}^{k_1-1} a_n \leq a_{k_1-1}^{k_1-1} < a_{k_1}^{k_1} \leq z^{k_1^2} < z^{N^2},$$

since $a_n$ is increasing and $k_1 < N$. Thus

$$\prod_{n=1}^{N-1} |a_n| < z^{N^2 + \left(\frac{M}{1-\beta}+\delta\right)^{-1} \left(\frac{M}{1-\beta}+1+\delta\right)^N}. \tag{28}$$

Since $\gamma$ was chosen as in Lemma 6.1, we have for each sufficiently large $z$ (and thereby $k_2$) that

$$\sum_{n=k_2}^{\infty} \frac{b_n}{a_n} \leq \frac{1}{a_{k_2}^\gamma}.$$

Combining this with inequality (25) and the fact that $N \leq k_2$, we find that

$$\sum_{n=k_2}^{\infty} \frac{b_n}{a_n} \leq \frac{1}{z^{\left(\frac{M}{1-\beta}+1+\delta\right)^{k_2}}} \leq \frac{1}{z^{\left(\frac{M}{1-\beta}+1+\delta\right)^N}} \tag{29}$$

when $z$ is sufficiently large. Since

$$1 - \beta = \frac{M}{\left(\frac{M}{1-\beta}\right)} > \frac{M}{\frac{M}{1-\beta}+\delta},$$

we may pick a fixed number $\zeta$ such that

$$\frac{M}{M/(1-\beta)+\delta} < \zeta < 1-\beta. \tag{30}$$

Now pick $0 < \Gamma < 1$ as in Lemma 6.3. Since $k_1$ is the largest number less than $k_2$ satisfying inequality (26) and $k_1 < N \leq k_2$, it follows that

$$\sum_{n=N}^{k_2-1} \frac{b_n}{a_n} \leq \frac{2^{\log_2^\Gamma a_N}}{a_N^{1-\beta}}$$

when $z$ (and thereby $N$) is sufficiently large. Applying inequalities (30) and (27), this yields that for all sufficiently large $z$,

$$\sum_{n=N}^{k_2-1} \frac{b_n}{a_n} < \frac{1}{a_N^\zeta} \leq \frac{1}{z^{\zeta\left(\frac{M}{1-\beta}+1+\delta\right)^N}}.$$

Combined with inequality (29), we conclude

$$\sum_{n=N}^{\infty} \frac{b_n}{a_n} < \frac{1}{z^{\zeta\left(\frac{M}{1-\beta}+1+\delta\right)^N}} + \frac{1}{z^{\zeta\left(\frac{M}{1-\beta}+1+\delta\right)^N}} < \frac{2}{z^{\zeta\left(\frac{M}{1-\beta}+1+\delta\right)^N}}, \tag{31}$$

when $z$ is sufficiently large.

Since $N$ is the smallest number strictly greater than $k_1$ satisfying inequality (27), and since $k_1$ cannot satisfy inequality (27), due to inequality (26), we get in particular that $N-1$ does not satisfy inequality (27). Thus,

$$2^{N^2 \log_2^c a_{N-1}} \leq 2^{N^2 \left(\frac{M}{1-\beta}+1+\delta\right)^{(N-1)c} \log_2^c z} < z^{N^2 \left(\frac{M}{1-\beta}+1+\delta\right)^{Nc}}. \tag{32}$$

When $z$ is sufficiently large, we obtain from inequalities (28), (31), and (32) that

$$Z_N = 2^{N^2 \log_2^c a_{N-1}} \left(\prod_{n=1}^{N-1} a_n^M\right) \sum_{n=N}^{k_2} \frac{b_n}{a_n}$$
$$< z^{N^2 \left(\frac{M}{1-\beta}+1+\delta\right)^{Nc} + MN^2 + \left(\frac{M}{M/(1-\beta)+\delta}-\zeta\right)\left(\frac{M}{1-\beta}+1+\delta\right)^N}.$$

To simplify notation, write $\zeta' = \zeta - \frac{M}{M/(1-\beta)+\delta}$, which is positive due to inequality (30). We then continue our calculation to find that

$$Z_N < z^{N^2 \left(\frac{M}{1-\beta}+1+\delta\right)^{Nc} + MN^2 - \zeta'\left(\frac{M}{1-\beta}+1+\delta\right)^N} < z^{-\frac{\zeta'}{2}\left(\frac{M}{1-\beta}+1+\delta\right)^N},$$

when $z$ is sufficiently large. As the right-hand side clearly tends to 0 as $z$ tends to infinity, we get the desired result.

**Case 2 (inequality (24) is not satisfied for any fixed $\delta > 0$)** This case is a bit more involved than the other one and will need to be split into two subcases, depending on whether $a_n \leq 2^n$ infinitely often. However, both cases will need an estimate of the expression

$$\frac{2^{n^2 \log_2^c a_{n-1} + \log_2^\Gamma a_n}}{(1+(n-1)^{-2})^{(1-\beta)\left(\frac{M}{1-\beta}+1\right)^n}}, \tag{33}$$

where $\Gamma \in (0, 1)$ is a fixed number to be chosen in each subcase. Set $\Gamma_0 = \max\{\Gamma, c\}$, and pick $\delta > 0$ so small that

$$\left(\frac{M}{1-\beta} + 1 + \delta\right)^{(1+\Gamma_0)/2} < \left(\frac{M}{1-\beta} + 1\right)^{(2+\Gamma_0)/3}. \tag{34}$$

When $n$ is sufficiently large, the case assumption will ensure that $a_n \leq 2^{\left(\frac{M}{1-\beta}+1+\delta\right)^n}$, and so it follows that

$$2^{n^2 \log_2^c a_{n-1} + \log_2^\Gamma |a_n|} \leq 2^{n^2 \left(\frac{M}{1-\beta}+1+\delta\right)^{cn} + \left(\frac{M}{1-\beta}+1+\delta\right)^{\Gamma n}} \leq 2^{\left(\frac{M}{1-\beta}+1+\delta\right)^{n(1+\Gamma_0)/2}}$$
$$< 2^{\left(\frac{M}{1-\beta}+1\right)^{n(2+\Gamma_0)/3}},$$

by applying inequality (34) in the last inequality. As for the denominator of expression. (33), the Taylor expansion of $\log_2(1 + x)$ implies that $\log_2(1 + n^{-2}) \geq n^{-5/2}$. Therefore,

$$\frac{2^{n^2 \log_2^c a_{n-1} + \log_2^\Gamma a_n}}{(1+n^{-2})^{(1-\beta)\left(\frac{M}{1-\beta}+1\right)^n}} < \frac{2^{\left(\frac{M}{1-\beta}+1\right)^{n(2+\Gamma_0)/3}}}{2^{n^{-5/2}(1-\beta)\left(\frac{M}{1-\beta}+1\right)^n}} \leq 2^{-n^{-3}\left(\frac{M}{1-\beta}+1\right)^n}, \tag{35}$$

for all sufficiently large $n$.

**Case 2.a** ($a_n \geq 2^n$ **for all but finitely many** $n$) By picking $\Gamma$ as in Lemma 6.2, we get for all sufficiently large $N \in \mathbb{N}$ that

$$\sum_{n=N}^{\infty} \left|\frac{b_n}{a_n}\right| \leq \frac{2^{\log_2^\Gamma |a_N|}}{|a_N|^{1-\beta}}.$$

At the same time, it follows from Lemma 6.5 and Lemma 6.4 with Eq. (16) that there are infinitely many $N \in \mathbb{N}$ such that

$$\prod_{n=1}^{N-1} |a_n|^M < \frac{\left(\max_{k \leq n < N} a_n^{\left(\frac{M}{1-\beta}+1\right)^{-n}}\right)^{(1-\beta)\left(\frac{M}{1-\beta}+1\right)^N}}{(1+N^{-2})^{(1-\beta)\left(\frac{M}{1-\beta}+1\right)^N}}$$
$$< \frac{|a_N|^{1-\beta}}{(1+N^{-2})^{(1-\beta)\left(\frac{M}{1-\beta}+1\right)^N}}.$$

Hence, for these infinitely many $N$,

$$Z_N = 2^{N^2 \log_2^c a_{N-1}} \left(\prod_{n=1}^{N-1} a_n^M\right) \sum_{n=N}^{k_2} \frac{b_n}{a_n} < \frac{2^{N^2 \log_2^c a_{N-1} + \log_2^\Gamma |a_N|}}{(1+N^{-2})^{(1-\beta)\left(\frac{M}{1-\beta}+1\right)^N}}.$$

From this and inequality (35), we obtain that for infinitely many $N$,

$$Z_N < 2^{-N^{-3}\left(\frac{M}{1-\beta}+1\right)^N},$$

and we are done.

**Case 2.b** ($a_n < 2^n$ **infinitely often**) Let $z > 0$ be sufficiently large, and pick $k_1, k_2, N \in \mathbb{N}$ as follows. Let $k_2$ be the smallest integer such that

$$a_{k_2}^{\left(\frac{M}{1-\beta}+1\right)^{-k_2}} > z, \tag{36}$$

and let $k_1$ be the largest integer such that $k_1 < k_2$ and

$$a_{k_1} < 2^{k_2}. \tag{37}$$

Due to the assumption that $a_n < 2^n$ infinitely often and the fact that $k_2$ is clearly unbounded, $k_1$ is also unbounded. Applying Lemma 6.4 with $k = k_1$ to Eq. (16), we pick $N$ to be the smallest integer such that $N > k_1$ and

$$a_N^{\left(\frac{M}{1-\beta}+1\right)^{-N}} > \left(1+N^{-2}\right) \max_{k_1 \leq n < N} a_n^{\left(\frac{M}{1-\beta}+1\right)^{-n}}. \tag{38}$$

Whenever $k_1 < n < N$, we then find by induction that

$$a_n^{\left(\frac{M}{1-\beta}+1\right)^{-n}} \leq \left(1+n^{-2}\right) \max_{k_1 \leq m < n} a_m^{\left(\frac{M}{1-\beta}+1\right)^{-m}} \leq \cdots$$

$$\leq \left(\prod_{m=k_1+1}^{n} \left(1+m^{-2}\right)\right) a_{k_1}^{\left(\frac{M}{1-\beta}+1\right)^{-k_1}}$$

$$\leq a_{k_1}^{\left(\frac{M}{1-\beta}+1\right)^{-k_1}} \prod_{m=1}^{\infty} \left(1+m^{-2}\right). \tag{39}$$

Since $\log(1+x) \leq x$, we find that

$$\prod_{m=1}^{\infty} \left(1+m^{-2}\right) = \exp\left(\sum_{m=1}^{\infty} \log\left(1+m^{-2}\right)\right) \leq \exp\left(\sum_{m=1}^{\infty} m^{-2}\right) = \exp\left(\frac{\pi^2}{6}\right) < 6.$$

Similarly, inequality (37) and the fact that $2^{n/2^n} \leq 1/2$ allow us to deduce

$$a_{k_1}^{\left(\frac{M}{1-\beta}+1\right)^{-k_1}} \leq 2^{k_1 2^{-k_1}} \leq \frac{1}{2}.$$

Recalling inequality (39), it follows that

$$a_n^{\left(\frac{M}{1-\beta}+1\right)^{-n}} < \frac{1}{2} \cdot 6 = 3, \tag{40}$$

for each $k_1 \leq n < N$, and so, due to Eq. (23) with $\delta = 0$, it follows that

$$\prod_{n=k_1}^{N-1} a_n^M < 3^{M \sum_{n=k_1}^{N-1} \left(\frac{M}{1-\beta}+1\right)^n} < 3^{\left(\frac{M}{1-\beta}+1\right)^N}$$

when $z$ (and thereby $N$) is sufficiently large. Using inequality (37) and the fact that $a_n$ is non-decreasing to estimate

$$\prod_{n=1}^{k_1-1} a_n^M \leq a_{k_1}^{Mk_1} \leq 2^{Mk_1^2} < 2^{MN^2}, \tag{41}$$

we may then conclude that

$$\prod_{n=1}^{N-1} a_n = \left(\prod_{n=1}^{k_1-1} a_n\right) \prod_{n=k_1}^{N-1} a_n < 2^{MN^2} 3^{\left(\frac{M}{1-\beta}+1\right)^N} \leq 4^{\left(\frac{M}{1-\beta}+1\right)^N}, \tag{42}$$

for all sufficiently large values of $z$. On the other hand, we might also estimate $\prod_{n=k_1}^{N-1} a_n^M$, using Lemma 6.5 and inequality (38) instead, which leads to

$$\prod_{n=k_1}^{N-1} a_n^M \leq \frac{a_N^{1-\beta}}{\left(1+N^{-2}\right)^{(1-\beta)\left(\frac{M}{1-\beta}+1\right)^N}},$$

so that we, by means of inequality (41), reach

$$\prod_{n=1}^{N-1} a_n^M \leq 2^{MN^2} \frac{a_N^{1-\beta}}{\left(1+N^{-2}\right)^{(1-\beta)\left(\frac{M}{1-\beta}+1\right)^N}}, \tag{43}$$

for all sufficiently large $z$ (and thereby $N$).

Let $0 < \gamma < 1$ be given as in Lemma 6.1. From this, inequality (36), and the fact that $N \leq k_2$ when $z$ is large, we then obtain

$$\sum_{n=k_2}^{\infty} \frac{b_n}{a_n} \leq a_{k_2}^{-\gamma} < z^{-\gamma\left(\frac{M}{1-\beta}+1\right)^{k_2}} \leq z^{-\gamma\left(\frac{M}{1-\beta}+1\right)^N}, \tag{44}$$

when $z$ is sufficiently large. Let similarly $0 < \Gamma < 1$ be given as in Lemma 6.3. Since $k_1$ is the largest number less than $k_2$ that satisfies inequality (37) and $k_1 < N \leq k_2$, we have

$$\sum_{n=N}^{k_2-1} \frac{b_n}{a_n} \leq \frac{2^{\log_2^\Gamma a_N}}{a_N^{1-\beta}}.$$

Together with inequality (44), this leads to

$$\sum_{n=N}^{\infty} \frac{b_n}{a_n} = \sum_{n=N}^{k_2-1} \frac{b_n}{a_n} + \sum_{n=k_2}^{\infty} \frac{b_n}{a_n} \leq \frac{2^{\log_2^\Gamma a_N}}{a_N^{1-\beta}} + z^{-\gamma\left(\frac{M}{1-\beta}+1\right)^N}. \tag{45}$$

Combining inequalities (42), (43), and (45), we obtain

$$\begin{aligned}
Z_N &= 2^{N^2 \log_2^c a_{N-1}} \left(\prod_{n=1}^{N-1} a_n^M\right) \sum_{n=N}^{k_2} \frac{b_n}{a_n} \\
&< 2^{N^2 \log_2^c a_{N-1}} \min\left\{4^{\left(\frac{M}{1-\beta}+1\right)^N}, 2^{MN^2} \frac{a_N^{1-\beta}}{\left(1+N^{-2}\right)^{(1-\beta)\left(\frac{M}{1-\beta}+1\right)^N}}\right\} \\
&\quad \cdot \left(\frac{2^{\log_2^\Gamma a_N}}{a_N^{1-\beta}} + z^{-\gamma\left(\frac{M}{1-\beta}+1\right)^N}\right) \\
&\leq \frac{2^{N^2 \log_2^c a_{N-1}+MN^2+\log_2^\Gamma a_N}}{\left(1+N^{-2}\right)^{(1-\beta)\left(\frac{M}{1-\beta}+1\right)^N}} + \frac{2^{N^2 \log_2^c a_{N-1}} 4^{\left(\frac{M}{1-\beta}+1\right)^N}}{z^{\gamma\left(\frac{M}{1-\beta}+1\right)^N}},
\end{aligned} \tag{46}$$

for all sufficiently large $z$. From inequality (35), it follows that

$$\frac{2^{N^2 \log_2^c a_{N-1}+MN^2+\log_2^\Gamma a_N}}{\left(1+N^{-2}\right)^{(1-\beta)\left(\frac{M}{1-\beta}+1\right)^N}} < \frac{2^{MN^2}}{2^{N-3\left(\frac{M}{1-\beta}+1\right)^N}},$$

which clearly tends to 0 as $z$ (and thereby $N$) tends to infinity. We are thus left to show that the remaining term of the right-hand side of inequality (46) also approaches 0 when $z$ grows large. Fortunately, this immediately follows from the calculation that

$$\frac{2^{N^2}\log_2^c a_{N-1} 4^{\left(\frac{M}{1-\beta}+1\right)^N}}{z^{\gamma\left(\frac{M}{1-\beta}+1\right)^N}} \leq \frac{5^{\left(\frac{M}{1-\beta}+1\right)^N}}{z^{\gamma\left(\frac{M}{1-\beta}+1\right)^N}} \leq z^{-\frac{\gamma}{2}\left(\frac{M}{1-\beta}+1\right)^N},$$

for all sufficiently large $z$. This completes the proof. □

## 7 Concluding remarks

In [12], the present author proves a variant of the main theorem of [6]. In particular, this implies that one may replace $\limsup_{n\to\infty} |a_n|^{\prod_{i=1}^{n-1}(d^i+d)^{-1}} = \infty$ with

$$\liminf_{n\to\infty} |a_n|^{\prod_{i=1}^{n-1}(d^i+d)^{-1}} < \limsup_{n\to\infty} |a_n|^{\prod_{i=1}^{n-1}(d^i+d)^{-1}} < \infty$$

in Theorem 1.3 and still get irrationality. Certainly, a corresponding result can be proven for Propositions 4.3 and 5.3 as well, leading to alternative versions of Theorems 1.4, 1.6, and 1.7, though we will not do that here.

We will now compare the main theorems to Theorem 1.2. If we set $d = 1$ in Theorems 1.6 or Theorem 1.7, we get the below corollary, which also appears if we assume $d = 1$ in the proof of Theorem 1.4.

**Corollary 7.1** *Let $\alpha, \delta, \varepsilon > 0$ be positive real numbers with $\alpha < 1$, and let $\beta \in [0, \frac{\varepsilon}{1+\varepsilon})$. Let $\{a_n\}_{n=1}^\infty$ and $\{b_n\}_{n=1}^\infty$ be sequences of positive integers so that*

$$n^{1+\varepsilon} \leq a_n \leq a_{n+1}, \quad \limsup_{n\to\infty} a_n^{\left(\frac{2+\delta}{1-\beta}+1\right)^{-n}} = \infty,$$

*and for all sufficiently large $n$,*

$$b_n \leq a_n^\beta 2^{\log_2^\alpha a_n}. \tag{47}$$

*Then the sequence $\{a_n/b_n\}_{n=1}^\infty$ is transcendental.*

This corollary also follows from Theorem 1.2 by replacing $\varepsilon$ with $\beta/(1-\beta) + \delta/3$ and putting $\gamma = \frac{2\beta+\delta}{1-\beta}$, while noting that $a_n^\beta 2^{\log_2^{\alpha_1} a_n} < a_n^{\varepsilon/(1+\varepsilon)}/2^{\log_2^{\alpha_2} a_n}$ for any fixed values of $\alpha_1, \alpha_2 \in (0, 1)$ and all sufficiently large $n$. On the other hand, Theorem 1.2 is slightly stronger than Corollary 7.1 since inequality (1) allows $\log |b_n|/\log |a_n|$ to approach $\varepsilon/(1+\varepsilon)$ as $n \to \infty$, which is prevented by inequality (47). Quite naturally, this raises the following question.

**Question 7.2** *Suppose we replaced $\beta$ by $\varepsilon/(1+\varepsilon)$, let $\alpha$ be sufficiently close to 1, and replaced the assumption $|b_n| \leq |a_n|^\beta 2^{\log_2^\alpha a_n}$ by $|b_n| < |a_n|^{\varepsilon/(1+\varepsilon)} 2^{-\log_2^\alpha |a_n|}$ for all sufficiently large n. Would Theorems 1.4, 1.6, and 1.7 then remain true?*

While the above comparison between Corollary 7.1 and Theorem 1.2 would suggest an affirmative answer, this question is not so easily answered. To see why, we start by taking a brief look at Hančl's proof of Theorem 1.2 as presented in [5]. Part of the proof follows an argument much similar to Lemma 3.3. The most significant difference is that Hančl takes advantage of the fact that $\gamma$ can always be replaced by a smaller value $\gamma' > 2\varepsilon$ without affecting whether any assumption of the theorem is satisfied. This means that he only ever has to consider what corresponds to case 1 from the proof of Lemma 3.3. In the proofs of Theorems 4.1 and 5.1, we can make the same trick by replacing $\delta$ with a smaller positive number, so the proofs of these theorems should be easily modified to allow the changes proposed by the above question. The problem arises with Propositions 4.3 and 5.3 – here none of the parameters present in the limsup conditions can be reduced without strengthening some other assumption, and so we cannot apply Hančl's trick and must deal with some version of case 2 from the proof of Lemma 3.3. When imposing the changes from Question 7.2, this becomes a much more difficult task, and there appears to be no immediate way of modifying the proof of Lemma 3.3 so that all of case 2 is covered.

Therefore, at least until Question 7.2 is answered in the affirmative for at least one of the theorems, we cannot both get the large values of $|b_n|$ that Theorem 1.2 suggests without making the limsup conditions more strict.

**Acknowledgements**
I thank the Independent Research Fund Denmark for funding the research, and I thank my supervisor Simon Kristensen for proofreading the paper.

**Funding** Open access funding provided by Aarhus Universitet. This work received funding from the Independent Research Fund Denmark (Grant number 1026-00081B).

**Data availability** Links to supporting data http://doi.org/10.48550/arXiv.2308.15302.

Received: 22 February 2024 Accepted: 23 June 2024 Published online: 28 July 2024

**References**
1. Erdős, P., Graham, R.L.: Old and new problems and results in combinatorial number theory. In: Monographies de L'Enseignement Mathématique [Monographs of L'Enseignement Mathématique], vol. 28. Université de Genève, L'Enseignement Mathématique, Geneva (1980)
2. Hančl, J.: Transcendental sequences. Math. Slov. **46**(2–3), 177–179 (1996)
3. Erdős, P.: Some problems and results on the irrationality of the sum of infinite series. J. Math. Sci. **10**, 1–7 (1975)
4. Kolouch, O., Novotný, L.: Diophantine approximations of infinite series and products. Commun. Math. (2016). https://doi.org/10.1515/cm-2016-0006
5. Hančl, J.: Two criteria for trascendental sequences. Mat. (Catania) **56**(1), 149–160 (2001)
6. Andersen, S.B., Kristensen, S.: Arithmetic properties of series of reciprocals of algebraic integers. Mon. Mat. **190**(4), 641–656 (2019). https://doi.org/10.1007/s00605-019-01326-1
7. Hančl, J., Nair, R., Šustek, J.: On the lebesgue measure of the expressible set of certain sequences. Indag. Math. **17**(4), 567–581 (2006). https://doi.org/10.1016/S0019-3577(06)81034-7
8. Schmidt, W.M.: Diophantine Approximation. Springer, Berlin (1980)
9. Waldschmidt, M.: Diophantine approximation on linear algebraic groups. In: Grundlehren der Mathematischen Wissenschaften, vol. 328. Springer, Berlin (2000)
10. Bugeaud, Y.: Approximation by algebraic numbers. In: Cambridge Tracts in Mathematics, vol. 160. Cambridge University Press, Cambridge (2004)
11. Stewart, I.: Galois Theory, 3rd edn. Chapman & Hall/CRC Mathematics, Boca Raton (2004)
12. Laursen, M.L.: Algebraic degree of series of reciprocal algebraic integers. Rocky Mt. J. Math. **53**(2), 517–529 (2023). https://doi.org/10.1216/rmj.2023.53.517